\documentclass[12pt]{amsart}
\usepackage{amssymb,amsmath,amstext,amsthm,amsfonts}
\usepackage[dvips]{graphicx}
\usepackage{a4wide}

\newcommand{\leb}{\operatorname{Leb}}
\newcommand{\dist}{\operatorname{dist}}

\begin{document}

\newcommand{\mcup}{\mbox{$\bigcup$}}
\newcommand{\mcap}{\mbox{$\bigcap$}}

\def \RR {{\mathbb R}}
\def \ZZ {{\mathbb Z}}
\def \NN {{\mathbb N}}
\def \PP {{\mathbb P}}
\def \TT {{\mathbb T}}
\def \II {{\mathbb I}}
\def \JJ {{\mathbb J}}

\def \ra {\rightarrow }
 \def \wh {\widehat }
  \def \leb {\mbox{Leb} }
 \def \un{\underline }

 \def \ov {\overline}
 \def \supp {\mbox{supp}\, }
 \def \wlim {\mbox{$w^*$-}\lim_{n\ra\infty}\, }
 \def \distp{\mbox{d}_\PP }

 \def \be {\beta } \def \vfi {\varphi
}  \def \si {\sigma }
\def \vare {\varepsilon }

 \def \cf {\mathcal{F}}
 \def \cm {\mathcal{M}}
 \def \cn {\mathcal{N}}
 \def \cq {\mathcal{Q}}
 \def \cp {\mathcal{P}}
 \def \cc {\mathcal{C}}
 \def \ch {\mathcal{H}}
  \def \cs {\mathcal{S}}

\newcommand{\dem}{\begin{proof}}
\newcommand{\cqd}{\end{proof}}

\newcommand{\qand}{\quad\text{and}\quad}

\newtheorem{theorem}{Theorem}
\newtheorem{corollary}{Corollary}

\newtheorem*{Maintheorem}{Main Theorem}

\newtheorem{maintheorem}{Theorem}
\renewcommand{\themaintheorem}{\Alph{maintheorem}}
\newcommand{\cmt}{\begin{maintheorem}}
\newcommand{\fmt}{\end{maintheorem}}

\newtheorem{maincorollary}[maintheorem]{Corollary}
\renewcommand{\themaintheorem}{\Alph{maintheorem}}
\newcommand{\cmc}{\begin{maincorollary}}
\newcommand{\fmc}{\end{maincorollary}}

\newtheorem{T}{Theorem}[section]
\newcommand{\ct}{\begin{T}}
\newcommand{\ft}{\end{T}}

\newtheorem{Corollary}[T]{Corollary}
\newcommand{\cco}{\begin{Corollary}}
\newcommand{\fco}{\end{Corollary}}

\newtheorem{Proposition}[T]{Proposition}
\newcommand{\cpr}{\begin{Proposition}}
\newcommand{\fpr}{\end{Proposition}}

\newtheorem{Lemma}[T]{Lemma}
\newcommand{\cle}{\begin{Lemma}}
\newcommand{\fle}{\end{Lemma}}

\newtheorem{Sublemma}[T]{Sublemma}
\newcommand{\csle}{\begin{Lemma}}
\newcommand{\fsle}{\end{Lemma}}

\newtheorem{Remark}[T]{Remark}
\newcommand{\cre}{\begin{Remark}}
\newcommand{\fre}{\end{Remark}}

\newtheorem{Definition}[T]{Definition}
\newcommand{\cd}{\begin{Definition}}
\newcommand{\fd}{\end{Definition}}

\title[Markov structures for
non-uniformly expanding systems]{Markov structures and decay of
correlations for non-uniformly expanding dynamical systems}

\author{Jos\'e F. Alves}
\address{Departamento de Matem\'atica Pura, Faculdade de Ci\^encias do Porto\\
Rua do Campo Alegre 687, 4169-007 Porto, Portugal}
\email{jfalves@fc.up.pt}
\urladdr{http://www.fc.up.pt/cmup/home/jfalves}

\author{Stefano Luzzatto}
\address{Mathematics Department, Imperial College\\
180 Queen's Gate, London SW7, UK}
\email{stefano.luzzatto@ic.ac.uk}
\urladdr{http://www.ma.ic.ac.uk/$\sim$luzzatto}

\author{Vilton  Pinheiro}
\address{Departamento de Matem\'atica, Universidade Federal da Bahia\\
Av. Ademar de Barros s/n, 40170-110 Salvador, Brazil.}
\email{viltonj@ufba.br}

\date{May 2002}

\thanks{Work carried out at the  Federal University of
Bahia,  University of Porto and Imperial College, London.
Partially supported by CMUP, PRODYN, SAPIENS and UFBA}

\setcounter{tocdepth}{1}

\maketitle

\begin{abstract}
We consider  non-uniformly expanding maps
 on compact Riemannian manifolds of arbitrary dimension, possibly
 having discontinuities and/or critical sets,
 and show that under some general conditions they admit an induced
 Markov tower structure for which the decay of the return time function
 can be controlled in terms of the time generic points need to achieve
 some uniform expanding behavior.  As a consequence we obtain some
 rates for the decay of correlations of those maps and conditions for the
 validity of the Central Limit Theorem.
\end{abstract}

\tableofcontents


\section{Introduction and statement of results}

The purpose of this paper is to study the geometrical structure
and statistical properties of piecewise smooth dynamical systems which
satisfy some asymptotic expansion properties almost everywhere.
We begin with a discussion of the statistical properties we are
interested in, and the precise statement of our assumptions and results
concerning these properties.  We then state our main result on the
existence of an induced Markov map and present our main application to
class of two-dimensional non-uniformly expanding \emph{Viana maps}.

\subsection{Statistical properties}

One of the most powerful ways of describing the dynamical features of
systems, specially those having a very complicated geometrical and topological
structure of individual orbits, is through invariant probability measures.
Any such measure can be decomposed into ergodic components and, by a
simple application of Birkhoff's Ergodic Theorem, almost every initial
condition in each ergodic component has the same statistical
distribution in space.  On such a component, a map \( f \) is said to
be \emph{mixing} if
\[
|\mu(f^{-n}(A) \cap B) - \mu(A)\mu(B)|\to 0,\quad\text{when
$n\to\infty$},
\]
for any measurable sets \( A, B \).
Standard counterexamples show that in general there is no specific
rate at which this loss of memory occurs: it is always possible to
choose sets \( A \) and \( B \) for which mixing is arbitrarily
slow. It is sometimes possible however, to define the
\emph{correlation function}
\[
\mathcal C_{n}(\varphi, \psi) = \left|\int (\varphi \circ f^{n})
\psi d\mu - \int \varphi d\mu \int \psi d\mu\right|,
\]
and to obtain specific rates of decay which depend only on the map
\( f \) (up to a multiplicative constant which is allowed to
depend on \( \varphi, \psi \)) as long as the \emph{observables}
\( \varphi, \psi \) belong to some appropriate functional space.
Notice that choosing these observables to be characteristic
functions this gives
exactly the original definition
of mixing.

The precise dynamical features which cause mixing, and in
particular the dynamical features which cause different rates of
decay of the correlation function, are still far from understood.
Exponential mixing for uniformly expanding and uniformly
hyperbolic systems has been known since the work of Sinai, Ruelle
and Bowen \cite{Sin68, Bow70, Bow75, BR} and may not seem surprising in
view of the fact that all quantities involved are exponential.
However the subtlety of the question is becoming more apparent in the
light of recent examples which satisfy asymptotic exponential
expansion estimates but only subexponential decay of correlations.
The simplest case is that of one-dimensional maps which are expanding
everywhere except at some fixed point \( p \) for which \( f'(p)=1 \).
In certain cases (essentially depending on the second derivative \(
f''(p) \)) there is an absolutely continuous mixing invariant measure
with positive Lyapunov exponent but strictly subexponential \cite{LSV99,
Hu99, Y2} (and in some cases even sub-polynomial \cite{Hol02}) decay of
correlations.  In this case the indifferent fixed point is
\emph{slowing down} the mixing process since nearby points are moving
away (and thus ``mixing'') at a slower, subexponential, rate rather
than the exponential rate at which they move away from other fixed or
periodic point.

A more subtle slowing down effect occurs in smooth one-dimensional
maps with critical points where the rate of mixing is essentially
determined by the rate of growth of the derivative along the
critical orbit \cite{BLS}.  Here, points close to the critical
point shadow its orbit for a certain amount of time slowing down
the mixing process like in the case of an indifferent fixed point
if the derivative growth along the critical orbit is
subexponential.  In this paper we identify for the first time a
general feature which plays an important role in determining the
rate of decay of correlation for the system.  This is the
\emph{degree of non-uniformity} of the expansivity which measures
how close the system is to being uniformly expanding by
quantifying the initial time one has to wait for typical points to
start behaving as though the system were uniformly expanding.  The
precise definition will be given below.

We also obtain conditions for the validity of the Central Limit
Theorem, which states that the probability of a given deviation of
the average values of an observable along an orbit from the
asymptotic average is essentially given by a Normal Distribution:
given a H\"older continuous function \( \phi \) which is not a
coboundary (\( \phi\neq \psi\circ f - \psi \) for any \( \psi \))
there exists \( \sigma>0 \) such that for every interval \(
J\subset \mathbb R \),
\[
\mu\left\{x\in X:
\frac{1}{\sqrt n}\sum_{j=0}^{n-1}\left(\phi(f^{j}(x))-\int\phi d\mu
\right)\in J \right\} \to
\frac{1}{\sigma \sqrt{2\pi} }\int_{J} e^{-t^{2}/ 2\sigma^{2}}dt.
\]
We present our results first of all in the case of local
diffeomorphisms and then in the case in which the map might contain
discontinuities and/or critical points.

\subsubsection*{Non-uniformly expanding local diffeomorphisms.}

Let $M$ be a compact Riemannian manifold  of dimension \( d\geq 1
\) and \( \leb \) a normalized Riemannian volume form on $M$ that
we call {\em Lebesgue measure}. Let \( f: M \to M \) be a \( C^{2}
\) local diffeomorphism and suppose that there exists a  constant
\( \lambda > 0 \)
 such that for Lebesgue almost all points \( x\in M \) the following
 \emph{non-uniform expansivity} condition is satisfied:
\begin{equation}\tag{\(*\)}
    \liminf_{n\to\infty}\frac{1}{n}\sum_{i=0}^{n-1}
    \log \|Df_{f^{i}(x)}^{-1}\|^{-1}\geq \lambda >0.
\end{equation}
  Notice that in the one-dimensional case this condition
reduces to
\[
 \liminf_{n\to\infty}\frac{1}{n}\sum_{i=0}^{n-1}
 \log |f'(f^{n}(x))|= \liminf_{n\to\infty} |(f^{n})'(x)| \geq \lambda >0.
\]
The formulation in the higher dimensional case is motivated by the fact
that we want to make an assumption about the average expansion \emph{in
every direction}.  Indeed for a linear map \( A: \mathbb R^{d}\to
\mathbb R^{d} \), the condition \( \|A\| > 1 \)
only provides
information about the existence of \emph{some} expanded direction,
whereas the condition \( \|A^{-1}\|^{-1} >1 \) (i.e. \( \log \|A^{-1}\|^{-1}> 0 \))
is exactly equivalent
to saying that \emph{every} direction is expanded by \( A \).
Condition \( (*) \) implies that the \emph{expansion time} function
\[
\mathcal E(x) =
\min\left\{N: \frac{1}{n}
\sum_{i=0}^{n-1} \log \|Df^{-1}_{f^{i}(x)}\|^{-1} \geq \lambda/2
\ \ \forall n\geq N\right\}.
\]
is defined and finite almost everywhere in \( M \). We think of
this as the \emph{waiting time} before the exponential derivative
growth kicks in.  Our results indicate that a main factor
influencing the rate of decay of correlation is rate of decay of
the tail of this function, i.e. the rate of decay of the measure
of the set of points which have not yet started expanding
uniformly by time \( n \). We remark that the choice of \(
\lambda/2 \) in the definition of \( \mathcal E \) is just for
convenience, any other positive constant smaller than \( \lambda
\) would work and would yield the same results.

\begin{theorem}\label{t:locdif}
    Let \( f: M\to M \) be a transitive \( C^{2} \) local diffeomorphism satisfying
    condition \( (*) \)  and suppose that there exists \( \gamma>1 \) such that
    \[
    \leb\big(\{\mathcal E(x)> n\}\,\big)\leq \mathcal O(n^{-\gamma}).
    \]
    Then there exists an absolutely continuous, \( f \)-invariant,
    probability measure \( \mu \). 
    Some finite power of \( f \) is
    mixing with respect to $\mu$ and the correlation function \( \mathcal C_{n} \) for H\"older
    continuous observable on \( M \) satisfies
    \[
    \mathcal C_{n} \leq \mathcal O(n^{-\gamma+1 }).
    \]
    Moreover, if \( \gamma > 2 \) then the Central Limit Theorem holds.
 \end{theorem}
 The existence of a measure \( \mu \) and the finiteness of the
 number of ergodic components of \( \mu \) was proved in \cite{ABV}.
 Our construction gives an alternative proof of the absolute
 continuity of \( \mu \) and allows us to obtain the estimates on the
 rate of Decay of Correlation and on the validity of the Central
 Limit Theorem.  We remark that the questions concerning existence and
 ergodicity of an absolutely continuous invariant measure are quite
 distinct from the questions of the statistical properties with respect to the
 measure.  Our results apply and are of interest even if an absolutely
 continuous, ergodic, \( f \)-invariant, probability measure \( \mu \)
 on \( M \) is already known to exist.  In fact, in this case condition (\( * \))
 admits a very natural formulation simply in terms of the average value of the
 (inverse of the) norm of the (inverse of the) derivative:
\begin{equation*}
    \int \log
\|Df^{-1}\|^{-1} d\mu > 0
\end{equation*}
Indeed Birkhoff's Ergodic Theorem then implies that the limit
\[
\lambda = \lim_{n\to\infty} \frac{1}{n}\sum_{i=0}^{n-1} \log \|Df^{-1}_{f^{i}(x)}\|^{-1}
=     \int \log
\|Df^{-1}\|^{-1} d\mu > 0
\]
exists for \( \mu \)-almost every \( x\in M \).  In particular
the expansion time function \( \mathcal E(x) \) is also defined and
finite almost everywhere and the conclusions of the Theorem hold under
the given conditions on the rate of decay of the measure of \(
\{\mathcal E(x) > n\} \).

\subsubsection*{Maps with critical points and discontinuities.}

We can generalize our results to the case in which \( f \) is a
local diffeomorphism outside a \emph{critical/singular} set \( \cs
\subset M\) satisfying the following geometrical non-degeneracy
conditions which essentially say that
 $f$ {\em behaves like a power of the distance}
 to \( \cs \):
 there are constants $B>1$ and $\be>0$ such that for every $x\in
 M\setminus\cs$
\begin{enumerate}
 \item[(S1)]
\quad $\displaystyle{\frac{1}{B}\dist(x,\cs)^{\be}\leq
\frac{\|Df(x)v\|}{\|v\|}\leq B\dist(x,\cs)^{-\be}}$ for all $v\in
T_x M$;
\end{enumerate}
Moreover the functions \(  \log\det Df \) and \( \log \|Df^{-1}\|
\) are \emph{locally Lipschitz} at points \( x\in M \setminus
\mathcal S \) with Lipschitz constant depending on $ \dist (x,
\mathcal S)\): for every $x,y\in M\setminus \cs$ with
$\dist(x,y)<\dist(x,\cs)/2$ we have
\begin{enumerate}
\item[(S2)] \quad $\displaystyle{\left|\log\|Df(x)^{-1}\|-
\log\|Df(y)^{-1}\|\:\right|\leq
B\frac{\dist(x,y)}{\dist(x,\cs)^{\be}}}$;
 \item[(S3)]
\quad $\displaystyle{\left|\log|\det Df(x)^{-1}|- \log|\det
Df(y)^{-1}|\:\right|\leq B\frac{\dist(x,y)}{\dist(x,\cs)^{\be}}}$;
 \end{enumerate}
In the
one-dimensional case,  points in \( \cs \) may be points for which
\( f'=0 \) and/or \( f \) is discontinuous and/or  \( f'=\pm
\infty \).  Notice that in this case conditions \( (S2) \) and \(
(S3) \) are equivalent and therefore the geometrical conditions
just say that the derivative which may tend to 0 or \( \infty \)
at \( \cs \) can do so at most at some uniform polynomial rate and
that the Lipschitz constant of the log of the derivative is
proportional to some power of the distance to the critical set.

We assume that orbits have {\em slow approximation} or {\em
subexponential recurrence} to the critical set in the following
sense. Let \( d_{\delta}(x,\cs) \) denote the \( \delta
\)-\emph{truncated} distance from \( x \) to \( \cs \) defined as
\( d_{\delta}(x,\cs) =  d(x,\cs) \) if \( d(x,\cs) \leq \delta\)
and \( d_{\delta}(x,\cs) =1 \) otherwise. Then, given any
$\epsilon>0$ there exists $\delta>0$ such that for Lebesgue almost
every $x\in M$
\begin{equation} \label{e.faraway1}\tag{\( ** \)}
    \limsup_{n\to+\infty}
\frac{1}{n} \sum_{j=0}^{n-1}-\log \dist_\delta(f^j(x),\cs)
\le\epsilon.
\end{equation}
Again this is an asymptotic statement and we have no a-priori
knowledge about how fast this limit is approached or with what degree
of uniformity for different points \( x \). Since some control of the
recurrence at finite times is important for our construction we introduce the
\emph{recurrence time} function
\[
\mathcal R(x) = \min\left\{N\ge 1: \frac{1}{n} \sum_{i=0}^{n-1}
-\log \dist_\delta(f^j(x),\cs) \leq 2\varepsilon, \ \ \forall
n\geq N\right\}.
\]
Condition \( (**) \) implies that the \emph{recurrence time}
function is defined and finite almost everywhere in \( M \).
Before we state our results in this case, it will be useful to
introduce for each $n\ge1$ the set
\[
    \Gamma_n=\{x: \mathcal E(x) > n \ \text{ or } \ \mathcal R(x) > n \}.
\]
This is the set of points which at time \( n \) have not yet achieved
either the uniform exponential growth or the uniform subexponential
recurrence given by conditions \( (*) \) and \( (**) \).

\begin{theorem}\label{t:crit}
    Let \( f:M\to M \) be a transitive \( C^{2} \) local diffeomorphism outside a
    critical/singular set \( \cs \) satisfying the non-degeneracy
    conditions stated above.  Suppose that \( f \) satisfies the
    non-uniform expansivity condition \( (*) \) and the slow approximation
    condition \((**)\) to the critical set and suppose that there exists
    \( \gamma>1 \) such that
    \[
    \leb(\Gamma_n) \leq \mathcal O(n^{-\gamma}).
    \]
    Then there exists an absolutely continuous, \( f \)-invariant,
    probability measure \( \mu \).
    Some finite power of \( f \) is
    mixing with respect to $\mu$ and the correlation function \( \mathcal C_{n} \) for H\"older
    continuous observables on \(M \) satisfies
    \[
    \mathcal C_{n} \leq \mathcal O(n^{-\gamma+1 }).
    \]
    Moreover, if \( \gamma > 2 \) then the Central Limit Theorem holds.
    \end{theorem}
Notice that the assumptions of Theorem \ref{t:crit} contain those
of Theorem \ref{t:locdif} as a special case where \( \cs =
\emptyset \).  We have stated the two results seperately because
the local diffeomorphism case is sufficiently interesting on its
own and to emphasize the fact that the recurrence condition only
applies to the case in which a critical and/or singular set
exists. Both Theorems extend to arbitrary dimension the results of
\cite{ALP1} in which similar results were obtained for
one-dimensional maps.

We remark also that even though condition (\(**\)) is not needed
in all its strength for the proof (it is sufficient that the
statement holds for some \( \varepsilon \) sufficiently small
depending on expansivity rate \( \lambda \) and on the constants
\( B, \beta \) in the non-degeneracy conditions for the critical
set),  it is nevertheless more natural than it might appear at
first sight.  For example,  if an ergodic, \( f \)-invariant,
absolutely continuous probability measure \( \mu \) is given, then
this condition just amounts to supposing that this invariant
measure does not give too much weight to neighbourhoods of \( \cs
\) in the sense that
\[
\int \left| \log \dist(x, \cs)\right|  d\mu < \infty.
\]
Indeed, as for the expansivity condition, this immediately implies
(\(**\)) by Birkhoff's
Ergodic Theorem.  Notice moreover that this integrability condition
is satisfied if the singular set \( \mathcal S \) and the Radon-Nykodim
derivative of \( \mu \) with resect to Lebesgue satisfy some mild
regularity conditions.

 \subsection{Markov structure}

Our strategy for proving the results stated above is to establish
the existence of a \emph{Markov tower} structure:  a ball \(
\Delta \subset M \) and a countable partition \( \mathcal P \) of
\( \Delta \) into topological balls with smooth boundaries with
the property that each element \( U \) of \( \mathcal P \) has an
associated return time \( R(U) \) so that \( f^{R(U)}(U) = \Delta
\) with some uniform bounds on the volume distortion between one
return and the next.  Moreover we set up a combinatorial and
probabilistic argument which allows us to obtain estimates for the
tail \( \leb(\{R>n\}) \) of the return time function in terms of
the tail of the expansivity and recurrence functions defined
above.

\newpage

 \begin{Maintheorem}
     \label{t:Markov towers}
 Let \( f:M\to M \) be a transitive \( C^{2} \) local diffeomorphism outside a
    critical/singular set \( \cs \) satisfying the non-degeneracy
    conditions stated above.  Suppose that \( f \) satisfies the
    non-uniform expansivity condition \( (*) \) and the slow approximation
    condition~\((**)\) to the critical set and suppose that there exists
      \( \gamma>0 \) such that
    \[
    \leb(\Gamma_n) \leq \mathcal O(n^{-\gamma}).
    \]
    Then there exists a ball  \( \Delta \subset M\setminus
    \mathcal S \), a countable partition \( \mathcal P \) (mod 0) of \(
    \Delta \) into topological balls \( U \) with smooth boundaries,
    and a return time function \( R: \Delta \to \mathbb N \) piecewise
    constant on elements of  \(\mathcal P \) satisfying the following
    properties:
    \begin{enumerate}
    \item {\em {Markov}:} for each \( U\in\mathcal P \) and
    \( R=R(U)
    \),
    \(
    f^{R}: U \to \Delta
    \)
    is a \( C^{2} \) diffeomorphism (and in particular a bijection).  Thus the
    induced map
    \[ F: \Delta \to \Delta \ \text{  given by } \
     F(x) = f^{R(x)}(x)
     \] is
    defined almost everywhere and satisfies the classical Markov property.

    \item {\emph{Uniform expansivity}:} There exists \( \hat\lambda > 1 \) such that
    for almost all \( x\in \Delta \) we have
    \[ \|DF(x)^{-1}\|^{-1}\geq \hat\lambda. \]
    In particular the separation time \( s(x,y) \) given by the maximum
    integer such that \( F^{i}(x) \) and \( F^{i}(y) \) belong to the same
    element of the partition \( \mathcal P \) for all \( i\leq s(x,y) \),
    is defined and finite for almost every pair of points \( x,y\in\Delta \).

    \item {\em {Bounded volume distortion}:}
    There exist a constant \( K >0 \) such that for any pair of points
    \( x,y\in \Delta \) with \( \infty > s(x,y) \geq 1\) we have
    \[
   \left|\frac{\det DF(x)}{\det DF(y)}-1\right| \leq
   K \hat\lambda^{-s(F(x),
    F(y))}.
    \]
    \item
    {\em {Polynomial decay of tail of return times}:} 
    \[ \leb(\{R>n\})\leq
    \mathcal O (n^{-\gamma}).  \]
\end{enumerate}
   \end{Maintheorem}

In the particular case that $\leb(\Gamma_n)$ decays faster than
any polynomial we obtain {\em super-polynomial decay} for the tail of
return times: $\leb(\{R>n\})\le\mathcal O(n^{-\gamma})$ for every
$\gamma>0$.

We remark that the significance of the existence of a Markov Tower
structures goes well beyond the consequences this has for the
statistical properties of the map.  It can be thought of as a partial
generalization, to the framework of non-uniformly expanding maps, of
the remarkable (and classical) Theorem of Bowen that any uniformly
hyperbolic compact invariant set for a \( C^{2} \) diffeomorphisms
admits a finite Markov partition (\cite{Bow70} see also
\cite{AW67,AW70,Sin68}).  Besides the intrinsic interest of such a
statement, this fact has been used innumerable times in relation to
all kinds of results concerning uniformly hyperbolic systems.
There has been some success in the direct generalization of this
result, for example to   systems with discontinuities
\cite{BSC91,KT92}.  However the constructions always give rise to
countable partitions and any conclusions about the invariant
measures and other statistical properties then depends on a
corresponding ergodic theory for countable subshifts which is much
less developed than the finite case, although some results exist,
see for example \cite{B99,M01, Sar}.

A significant break-through was achieved recently by Young in
\cite{Y1,Y2} where the general problem of proving the existence of
Markov partitions was essentially reformulated in terms of proving
the existence of \emph{Markov Towers} or \emph{induced Markov
maps} as defined above.  One important advantage of these
structures is that statistical information about the system is
deduced from statistical information about the return times and
not encoded in some kind of transition matrix which would in
general be very hard to determine.  Moreover the actual
construction of Markov Towers has at least two significant
advantages.  Firstly, one can choose conveniently some small
region of the dynamical phase space, instead of having to
construct a partition of the entire space, and use
\emph{approximate} information about the remaining part of the
space to construct a return map.  Secondly, one does not need  a
single iterate of the map to have special (Markov) properties, but
is allowed to wait a certain amount of time until this property is
obtained.  Most importantly, only some approximate (statistical)
information is required concerning the length of this waiting
time.

\subsection{Viana maps}
An important class of nonuniform expanding dynamical
    systems (with critical sets) in dimension greater than one was
    introduced by Viana in \cite{V}. This has served as a model
  for some relevant results on the ergodic
    properties of
    non-uniformly expanding maps in higher dimensions; see
    \cite{Al,AA,ABV,AV}.

 This class of  maps
 can be described as
follows. Let $a_0\in(1,2)$ be such that the critical point $x=0$
is pre-periodic for the quadratic map $Q(x)=a_0-x^2$. Let
$S^1=\RR/\ZZ$ and $b:S^1\rightarrow \RR$ be a Morse function, for
instance, $b(s)=\sin(2\pi s)$. For fixed small $\alpha>0$,
consider the map
 \[ \begin{array}{rccc} \hat f: & S^1\times\RR
&\longrightarrow & S^1\times \RR\\
 & (s, x) &\longmapsto & \big(\hat g(s),\hat q(s,x)\big)
\end{array}
 \]
 where  $\hat q(s,x)=a(s)-x^2$ with
$a(s)=a_0+\alpha b(s)$, and $\hat g$ is the uniformly expanding
map of the circle defined by $\hat{g}(s)=ds$ (mod $\ZZ$) for some
large integer $d$. In fact, $d$ was chosen greater or equal to 16
in \cite{V}, but recent results in \cite{Bu} showed that some
estimates in \cite{V} can be  improved and $d=2$ is enough. It is
easy to check that for $\alpha>0$ small enough there is an
interval $I\subset (-2,2)$ for which $\hat f(S^1\times I)$ is
contained in the interior of $S^1\times I$. Thus, any map $f$
sufficiently close to $\hat f$ in the $C^0$ topology has
$S^1\times I$ as a forward invariant region. We consider from here
on these maps restricted to $S^1\times I$.

Taking into account the expression of $\hat f$ it is not difficult
to check that it behaves like a power of the distance close to the
critical set $\{x=0\}$. Moreover, there is a small neighbourhood
$\cn$ of $\hat f$ in the $C^3$ topology of maps from $S^1\times I$
into itself, such that any  $f\in\cn$ also behaves like a power of
the distance close to its critical set, which is close to
$\{x=0\}$.
 The
most important results for $f\in \cn$ are summarized below:
\begin{enumerate}
\item $f$ is non-uniformly expanding and its orbits have slow
approximation to the critical set \cite{V,AA}; \item there are
constants $C,c>0$ such that $\leb(\Gamma_n)\le Ce^{-c\sqrt n}$ for
every $n\ge1$ \cite{V,AA}; \item $f$ is topologically mixing and
has a unique ergodic absolutely continuous
 invariant (thus SRB) measure \cite{Al,AV};
\item the density of the SRB measure varies continuously in the $L^1$ norm with $f$
\cite{AV};
\item $f$ is stochastically stable \cite{AA}.
\end{enumerate}

The decay of correlations for Viana maps remained unknown for
several years and one of the initial purposes of the present work
was trying to find rates of mixing for these maps. Fortunately
this led to results which hold in greater generality. As a
consequence of our theorems, we obtain the following result:

\begin{theorem}
    Any $f\in\cn$ has super-polynomial decay of correlations and Central Limit Theorem
    holds for $f$.
    \end{theorem}

\subsection{Remarks}



Before starting the proof of the Main Theorem we discuss our basic
strategy and the main technical issues involved in the construction.

\subsubsection*{Strategy.}
We start by choosing essentially arbitrarily a point \( p \) with
dense pre-images and some sufficiently small ball \( \Delta_{0} \)
around this point.  This will be the domain of definition of our
induced map. We then attempt to implement the naive strategy of
iterating \( \Delta_{0} \) until we find some good return iterate
\( n_{0} \) such that \( f^{n_{0}}(\Delta_{0}) \) completely
covers \( \Delta_{0} \) and some bounded distortion property is
satisfied.  There exists then some topological ball  \( U \subset
\Delta_{0} \) such that \( f^{n_{0}}(U) = \Delta_{0} \).  This
ball is then by definition an element of the final partition of \(
\Delta_{0} \) for the induced Markov map and has an associated
return time \( n_{0} \).  We then continue iterating the
complement \( \Delta_{0}\setminus U \) until more good returns
occur.  Most of the paper is dedicated to showing that this
strategy can indeed be implemented in a successful way, yielding a
partition (mod 0) of \( \Delta_{0} \) into piecewise disjoint
subsets, and an associated return time function which is Lebesgue
integrable.  The construction also yields substantial information
about the tail of the return time function, i.e. the decay of the
measure of the set of points whose return time is larger than \( n
\).  Indeed the main motivation for this paper is to show that the
rate of decay of this tail  is closely related to the rate at
which the derivative along orbits approaches the asymptotic
expansion rate.

\subsubsection*{Technical issues.}
There are two main technical difficulties, distinct but related to
each other, in carrying out the plan suggested above.  The first
has to do with the geometry of the returns to \( \Delta_{0} \),
and in particular of the geometry of the set of points which does
not return at a given time.  Such a set can be visualized as a
ball \( \Delta_{0}\) containing an increasing number of smaller
topological balls corresponding to the elements of the final partition
which have return times smaller than \( n \).  The exact location and
shape of these smaller disks is quite difficult to control, as is the
location and shape of their images at time \( n \).  Therefore some
care is required, as well as the introduction of some auxiliary
partitions and waiting times, to make sure that the set of points
returning at time \( n \) is disjoint from the set of points which
have already returned at some earlier time.  These geometrical issues
are essentially related to the higher dimensional nature of the
dynamics and arise also in the uniformly expanding case.  This case
has been treated in \cite{Y1} and we follow essentially the same
strategy and notation here.  We still give all the details, for
completeness and to make sure that any further problems associated to
the non-uniformity of the expansion are dealt with as well.

The second technical problem, on the other hand, is precisely due
to the strictly non-uniform nature of expansion in our situation.
The process of defining the set \( \{R\leq n\} \) of points which
have an associated return time less than or equal to \( n \), as a
union of disjoint topological disks in \( \Delta_{0} \),  gives
rise to very ``small'' regions in the complement \( \{R>n\} \),
i.e. regions which are squeezed into strange shapes by the
geometry of the previous returns. It is important to control the
extent to which this can happen and to show that even these small
regions eventually grow large enough so that they can cover \(
\Delta_{0} \) and thus contain an element of the partition for the
induced map.  In the uniformly hyperbolic case, once the suitable
definitions and notation have been introduced, a relatively
straightforward calculation shows this to be the case and shows
that in fact this growth of small regions to uniformly large scale
occurs uniformly exponentially fast.  In our context we only have
much more abstract information about the eventual expansion at
almost every point and therefore this part of the argument is more
subtle.

We shall use the idea of \emph{Hyperbolic Times} to show that our
assumptions imply that almost every point has a basis of
(arbitrarily small) neighbourhoods which at some time are mapped
to uniformly large scale with bounded distortion.  It follows that
the speed at which this large scale is achieved is not uniform but
rather depends on the distribution of hyperbolic times associated
to points in the regions in question, which can be arbitrarily
large. We conclude that the final return time function for the
Markov induced map is related to the statistics of hyperbolic
times.  Since hyperbolic times are naturally related to to speed
at which some uniform expansion estimates begin to hold, this
yields our desired conclusions.

One of the key issues we have to address is the relation between
the statistics of hyperbolic times, the spatial distribution of
points having hyperbolic time at some given time, and the
geometrical structure of sets arising from the construction of the
partition described above; see Corollary~\ref{c:hyperbolic3},
equation (\ref{e.media}) and Proposition~\ref{p.final}. We are
able to implement a partially successful strategy in this respect:
in the polynomial case we establish an essentially optimal link
between the rate of decay of the expansion/recurrence function and
the rate of decay of correlations. The nature of the argument does
not immediately extend to the exponential case.


\subsubsection*{Overview of the paper.}
The paper is completely dedicated to the proof of the Main Theorem on
the existence of the Markov tower and the associated tail estimates.  By recent
results of Young \cite{Y1,Y2} the rate of decay of the tail of the
return time function in this framework has direct implications for the
rates of decay of correlations and the Central Limit Theorem and
therefore Theorems \ref{t:locdif} and \ref{t:crit} follow by
an application of her results.

In section \ref{growing times} we give several estimates related
to the time it takes for small domains to grow to some fixed size
while preserving some bounded distortion properties.  In section
\ref{induced map} we give the precise algorithm for constructing
the Markov Tower and describe the associated combinatorial
information. The final three sections \ref{s.check},
\ref{transitional estimates} and \ref{s.decay} are dedicated to
proving that this algorithm effectively results in a countable
partition (mod 0) with the required properties.

\section{Growing to large scale}
\label{growing times}

In this section we give the basic growth estimates on which the
algorithm for the construction of the Markov Tower is based.
First of all we define the notion of Hyperbolic Time and show that
almost all points have an infinite basis of neighbourhoods which grow
to some fixed size with bounded distortion for some
corresponding infinite sequence of hyperbolic times.  The set
of hyperbolic times depends on the point and the first hyperbolic time
for a given point can be arbitrarily large in general, although we do
have some degree of control since it is related to the values of the
expansivity and recurrence functions \( \mathcal E \) and \( \mathcal
R \) at that point.  Next we prove a useful and non-obvious
consequence of our assumptions, namely that if we fix some \(
\varepsilon > 0 \) then there exist some \( N_{\varepsilon} \)
depending only on \( \varepsilon \) such that any ball of radius \(
\varepsilon \) has some subset which grows to a fixed size with
bounded distortion within \( N_{\varepsilon} \) iterates.  Finally we
show that our ``base'' \( \Delta_{0} \) can be chosen in such a way
that any other sufficiently large ball contains a subset which is
mapped bijectively to \( \Delta_{0} \) with bounded distortion and
within some fixed number of iterates.  A combination of these estimate
will play a crucial role in obtaining control of the tail of the
return times to \( \Delta_{0} \).

\subsubsection*{Hyperbolic times: growing to uniform scale in variable time.}

Let $B>1$ and $\beta>0$ be as in the hypotheses (S1)-(S4).
 In what follows $b$ is any fixed constant satisfying $0 < b <
\min\{1/2,1/(4\beta)\}$.
 Given $\sigma<1$ and
   $\delta>0$, we say that $n$
   is a {\em $(\sigma,\delta)$-hyperbolic time} for
a point $x\in M$ if  for all $1\le k \le n$, $$
\prod_{j=n-k}^{n-1}\|Df(f^j(x))^{-1}\| \le \sigma^k \qand
\dist_\delta(f^{n-k}(x), \cs)\ge \sigma^{b k}.
 $$
 For each \( n\geq 1
\) we define
\[
H_{n}=H_{n}(\sigma, \delta)=\{x\in M: n \text{ is a \(
(\sigma,\delta) \)-hyperbolic time for } x\}.
\]
We give two well-established results which show that \emph{i}) if
\( n \) is a hyperbolic time for \( x \), the map \( f^{n} \) is a
diffeomorphism with uniformly bounded volume distortion on a
neighborhood of \( x \) which is mapped to a disk of uniform
 radius; \emph{ii}) almost every
point has lots of hyperbolic times. We say that $f^n$  has
\emph{volume distortion} bounded by $ D$ on a set $V$ if, for
every $x,y\in V$,
$$ \frac{1}{  D} \le \frac{|\det Df^n
(x)|}{|\det Df^n (y)|}\le D \,.  $$

\cle \label{l.hyperbolic1} Given $\sigma<1$ and $\delta>0$, there
exist $\delta_1, D_{1},\kappa   >0$, depending only on
$\sigma,\delta$ and on the map $f$, such that for any $x\in M$ and
\( n\geq 1 \) a $(\sigma,\delta)$-hyperbolic time for \( x \),
there exists a neighborhood \( V_n(x) \) of \( x \) with the
following properties:
\begin{enumerate}
\item $f^{n}$ maps $V_n(x)$ diffeomorphically onto the ball
$B(f^{n}(x), \delta_1 )$;
\item for $1\le k <n$ and $y,
z\in V_n(x)$, $ \dist(f^{n-k}(y),f^{n-k}(z)) \le
\sigma^{k/2}\dist(f^{n}(y),f^{n}(z))$;
\item $f^n$ has volume distortion bounded by $ D_1$ on
$V_n(x)$;
\item 
\(V_n(x) \subset B(x, \kappa^{-n}) \).
\end{enumerate}
\end{Lemma}
\dem For the proofs of items 1, 2, 3 see Lemma 5.2 and Corollary
5.3 in \cite{ABV}.  Item~4  is an immediate consequence of  item
2. \cqd

We shall often refer to the sets \( V_n(x) \) as \emph{hyperbolic
pre-balls} and to their images \( f^{n}(V_n(x)) \) as
\emph{hyperbolic balls}.  Notice that the latter are indeed balls
of radius \( \delta_1 \).

\begin{Lemma}\label{l:hyperbolic2}
    There exists \( \theta>0 \) and $\delta>0$ depending only on $f$ and $\lambda$
     such that for Lebesgue almost every \( x
    \in M \) and \( n \geq \mathcal E(x) \) there exist $(\sigma,\delta)$-hyperbolic times
    $1 \le n_1 < \cdots < n_l \le n$ for \( x \) with $l\ge\theta
    n$.
    \end{Lemma}
\dem See Lemma 5.4 of \cite{ABV}. Let us remark  for the sake of
completeness that the proof of the lemma gives
$\sigma=e^{-\lambda/4}$. \cqd

\begin{Corollary}\label{c:hyperbolic3}
For every $n\ge1$  and every $A\subset M\setminus\Gamma_n$ with
positive Lebesgue measure we have
    \[
   \frac{1}{n}\sum_{j=1}^{n}\frac{\leb(A\cap H_{j})}{\leb(A)} \geq \theta.
   \]
\end{Corollary}
\begin{proof}
Take $n\ge 1$ and $A\subset M\setminus\Gamma_n$ with positive
Lebesgue measure. Observe that  $n\ge \mathcal E(x)$ for all $x\in
A$, by definition of $\mathcal E(x)$. Let $\xi_n$ be the measure
in $\{1,\dots,n\}$ defined by $\xi_n(J)=\# J/n$, for each subset
$J$. Then, using Fubini's theorem
\begin{eqnarray*}
\frac{1}{n} \sum_{j=1}^{n}\leb(A\cap H_j)
& =& \int \left(\int_A \chi(x,i)\,d\leb(x)\right)d\xi_n(i) \\
& = & \int_A \left(\int \chi(x,i)\,d\xi_n(i)\right)d\leb(x),
\end{eqnarray*}
where $\chi(x,i)=1$ if $x\in H_i$ and $\chi(x,i)=0$ otherwise.
Now, Lemma~\ref{l:hyperbolic2} means that the integral with
respect to $d\xi_n$ is larger than $\theta>0$. So, the last
expression above is bounded from below by $\theta \leb(A)$.
\end{proof}

\subsubsection*{Growing to uniform scale in uniform time.}
Now we show a simple (albeit slightly counterintuitive) fact that
any \( \vare \) ball  has a subset which grows to fixed size
within some uniformly bounded maximum number of iterates.

\begin{Lemma}\label{le:grow1}
    For each \( \varepsilon > 0 \) there exists \(  N_{\varepsilon}
    > 0 \)    such that any ball \( B
    \subset M \) of radius \( \varepsilon >0\) contains a hyperbolic pre-ball \(
    V_{n}\subset B \) with \( n\leq  N_{\varepsilon} \).
\end{Lemma}

\begin{proof}
    Given \( \varepsilon > 0 \) and a ball $B(z,\vare)$, choose \(  N'_{\varepsilon} \)
    large enough so that any hyperbolic pre-ball \( V_{n} \) associated to a
    hyperbolic time \( n\geq  N'_{\varepsilon} \) will be
    contained in a ball of radius \( \varepsilon/10 \) (\(
    N'_{\varepsilon} \sim \kappa^{-1}\log (10\varepsilon^{-1}) \)).  Now
    notice that each point has an infinite number of hyperbolic times and
    therefore we have that
    \[
    \leb\left(M\setminus \bigcup_{j= N'_{\varepsilon}}^{n}H_{j}\right) \to 0
    \quad\text{ as } n\to\infty.
    \]
    Therefore it is possible to choose
    \[
     N_{\varepsilon} = \min\left\{n\geq  N'_{\varepsilon}:
    \leb\left(M\setminus \bigcup_{j=
    N'_{\varepsilon}}^{n}H_{j}\right)\right\}
    \leq \varepsilon^{d}/10
    \]
    where \( d \) is the dimension of \( M \). This ensures that there is
    a point \( \hat x \in B(z, \varepsilon/2) \) with a hyperbolic time
    \( n\leq  N_{\varepsilon} \) and associated hyperbolic pre-ball \(
    V_{n}(x)\subset B(z, \varepsilon) \).
\end{proof}

\subsubsection*{Returning to a given domain.}

Now we derive an useful consequence of the transitivity of $f$.
Given $\delta>0$, we say that a subset $A$ of $M$ is {\em
$\delta$-dense} if any point in $M$ is at a distance smaller than
$\delta$ from $A$.

\cle Given $\delta>0$ there is $p\in M$ and $N_0\in\NN$ such that
$\bigcup_{j=0}^{N_0}f^{-j}(\{p\})$ is $\delta$~dense in $M$ and
disjoint from $\cs$. \fle

\dem Observe that the properties of $f$ imply that the images and
preimages of sets with zero Lebesgue measure still have zero
Lebesgue measure. Hence, the set
 \begin{equation}\label{e.set}\nonumber
 \mathcal B=\bigcup_{n\ge 0}f^{-n}\left(\bigcup_{m\ge
 0}f^{-m}(\cs)\right)
 \end{equation}
 has Lebesgue measure equal to zero.
 On the other hand, since $f$ is transitive, we have by \cite{ABV}
that there is a unique SRB measure for $\mu$, which is an ergodic
and absolutely continuous with respect to Lebesgue measure, and
whose support is the whole manifold $M$. Moreover, the ergodicity
of $\mu$ implies that $\mu$ almost every point in $M$ has a dense
orbit. Since $\mu$ is absolutely continuous with respect to
$\leb$, then there is a positive Lebesgue measure subset of points
in $M$ with dense orbit. Thus there must be some point $q\in
M\setminus \mathcal B$ with dense orbit. Take $N_0\in\NN$ for
which $q, f(q), \dots, f^{N_0}(q)$ is $\delta$-dense. The point
$p=f^{N_0}(q)$ satisfies the conclusions of the lemma.
 \cqd

We fix once and for all  \( p \in M \) and $N_0\in\NN$  for which
$\bigcup_{j=0}^{N_0}f^{-j}(\{p\})$ is $\delta_1/3$ dense in \( M
\) and disjoint from the critical set \( \cs \). Recall that
$\delta_1>0$ is the radius of hyperbolic pre-balls given by
Lemma~\ref{l.hyperbolic1}. Take constants $\vare>0$ and
$\delta_0>0$ so that
 $$
 \sqrt\delta_0 \ll \delta_1/2\qand 0< \varepsilon \ll \delta_0
 .
 $$

\begin{Lemma}\label{le:grow2}
There exist a constant \(   D_{0} 
>0\) depending only on \( f, \sigma, \delta_1 \) and the point~\(
p \), such that for any ball \( B\subset M \) of radius \(
\delta_1 \) (in particular for any hyperbolic ball), there exist
an open set \( V\subset B \) and an integer \( 0\leq m \leq
 N_{0} \) for which
   \begin{enumerate}
       \item \( f^{m} \) maps \( V \) diffeomorphically onto \(  B(p,
       2\sqrt\delta_{0}) \);
       \item \( f^{m}| V \) has volume distortion bounded by \(  D_{0}
       \). 
  \end{enumerate}
\end{Lemma}

\begin{proof}
Since $\bigcup_{j=0}^{N_0}f^{-j}(\{p\})$ is \( \delta_1/3 \) dense
in \( M \) and disjoint from $\cs$, choosing \( \delta_{0}>0 \)
sufficiently small we have that each connected component of the
preimages of \(B(p, 2\sqrt\delta_{0}) \) up to time \( N_{0} \)
are bounded away from the critical set \( \cs \) and are contained
in a ball of radius \( \delta_1/3 \).

This immediately implies that any ball \( B \subset M \) of radius
\( \delta_1 \) contains a preimage \( V \) of \( B(p,
2\sqrt\delta_{0}) \) which is mapped diffeomorphically onto \(
B(p, 2\sqrt\delta_{0}) \) in at most $ N_{0}$ iterates. Moreover,
since the number of iterations and the distance to the critical
region are uniformly bounded, the volume distortion is uniformly
bounded.
\end{proof}

\cre\label{r.c0} It will be useful to emphasize that  \(
\delta_{0} \) and $N_0$ have been chosen in such a way that  all
the connected component  of the preimages of \(B(p,
2\sqrt\delta_{0}) \) up to time \( N_{0} \) satisfy the
conclusions of the lemma. In particular, they are uniformly
bounded away from the critical set \( \cs \), and so there is some
constant $C_0>1$ depending only on $f$ and $\delta_1$ such that
 $$
 \frac1{C_0}\le \|Df^m(x)\|, \|(Df^m(x))^{-1}\|
 \le C_0
 $$
 for all $1\le m\le N_0$ and $x$ belonging to an $m$-preimage of $B(p,
 2\sqrt\delta_{0})$.
 \fre


\section{The partitioning algorithm}
\label{induced map}

We now describe the construction of the (mod 0) partition of \(
\Delta_{0}=B(p,\delta_{0}) \).  The basic intuition is that we
wait for some iterate \( f^{k}(\Delta_{0}) \) to cover \(
\Delta_{0} \) completely, and then define the subset \( U \subset
\Delta_{0} \) such that \( f^{k}: U \to \Delta_{0} \) is a
diffeomorphism, as an element of the partition with return time \(
k \).  We then continue to iterate the complement \(
\Delta_{0}\setminus U \) until this complement covers again \(
\Delta_{0} \) and repeat the same procedure to define more
elements of the final partition with higher return times.  Using
the fact that small regions eventually become large due to the
expansivity condition (and the lemmas given above), it follows
that this process can be continued and that Lebesgue almost every
point eventually belongs to some element of the partition and that
the return time function depends on the time that it takes small
regions to become large on average and this turns out to depend
precisely on the tail of the expansivity condition function.

The formalization of this argument requires several technical
constructions which we explain below. The construction is
inductive and we give precisely the general step of the induction
below.  For the sake of a better visualization of the process, and to
motivate the definitions, we start with the first step.

\subsubsection*{First step of the induction.}
  First of all we introduce
 neighborhoods of \( p \)
$$\Delta^{0}_{0}=\Delta_0=B(p,\delta_{0}),\quad
\Delta^{1}_{0}=B(p,2\delta_{0}),\quad
\Delta^{2}_{0}=B(p,\sqrt\delta_{0})\qand
\Delta^{3}_{0}=B(p,2\sqrt\delta_{0}).$$ For $0<\sigma<1$ given by
Lemma~\ref{l.hyperbolic1}, let
\[
I_{k}=\left\{x\in\Delta^{1}_{0}\: : \:\delta_{0}(1+\sigma^{k/2}) <
\dist(x,p) < \delta_{0}(1+\sigma^{(k-1)/2})\right\},\quad  k\ge 1,
\]
be a partition (mod 0) into countably many rings of \(
\Delta_0^{1}\setminus \Delta_0 \). Take \( R_{0} \) 
some large integer to be determined below; we ignore any dynamics
occurring up to time \( R_{0} \). Let \( k\geq R_{0}+1 \) be the
first time that \( \Delta_0\cap H_{k}\neq\emptyset \). For \( j <
k \) we define formally the objects \( \Delta_{j}, A_{j},
A_{j}^{\varepsilon} \) whose meaning will become clear in the next
paragraph, by \( A_{j}=A_{j}^{\varepsilon}=\Delta_{j}=\Delta_{0}
\). Let $(U_{k,j}^3)_j$ be the connected components of
 $
 f^{-k}(\Delta_0^3)\cap A_{k-1}^\vare
 $
 contained in hyperbolic pre-balls $V_{k-m}$ with $k- N_0\le
 m\le k$ which are mapped diffeomorphically onto $\Delta_0^3$ by $f^k$.
 Now let
\[
U_{k,j}^i=U_{k,j}^3\cap f^{-k}\Delta_0^i,\quad i=0,1,2
\]
and set $R(x)=k$ for $x\in U_{k,j}^0$. Now  take
\[
\Delta_{k}=\Delta_{k-1}\setminus \{R=k\}.
\]
We define also a function \( t_{k}:\Delta_{k}\to \mathbb N \) by
\begin{equation*}
    t_{k}(x) =
    \begin{cases}
    s & \text{ if }  x\in U_{k,j}^{1} \text{ and }
    f^{k}(x) \in I_{s} \text{  for some $j$;} \\
    0 & \text{ otherwise}.
    \end{cases}
\end{equation*}
Finally let
\[
A_{k}= \{x\in\Delta_{k}: t_{k}(x) = 0\}, \quad B_{k}=
\{x\in\Delta_{k}: t_{k}(x) > 0\}
\]
and
 $$
 A_{k}^\vare= \{x\in\Delta_{k}:
 \dist(f^{k+1}(x),f^{k+1}(A_k))<\vare\}.
 $$

\subsubsection*{General step of the induction.}

The general inductive step of the construction now follows by
repeating the arguments above with minor modifications.  More
precisely we assume that sets $\Delta_{i}$, $A_{i}$, $A_{i}^\vare$
 $B_{i}$, $\{R=i\}$ and functions \( t_{i}:
\Delta_{i}\to\mathbb N \) are defined for all \( i\leq n-1 \). For
\( i\leq R_{0} \) we just let \( A_{i}=A_{i}^\vare=\Delta_{i}=
\Delta_{0}\), \( B_{i}=\{R=i\}=\emptyset \) and \( t_{i}\equiv 0
\). Now let $(U_{n,j}^3)_j$ be the connected components of
 $
 f^{-n}(\Delta_0)\cap A_{n-1}^\vare
 $
 contained in hyperbolic pre-balls $V_m$, with $n-  N_0\le
 m\le n$, which are mapped onto $\Delta_0^3$ by $f^n$.
 Take
\[
U_{n,j}^i=U_{n,j}^3\cap f^{-n}\Delta_0^i,\quad i=0,1,2,
\]
and set $R(x)=n$ for $x\in U_{n,j}^0$. Take also
\[
\Delta_{n}=\Delta_{n-1}\setminus \{R=n\}.
\]
The definition of the function \( t_{n}:\Delta_{n}\to \mathbb N \)
is slightly different in the general case.
\begin{equation*}
    t_{n}(x) =
    \begin{cases}
    x & \text{ if } x\in U_{n,j}^{1}\setminus U_{n,j}^{0} \text{ and }
    f^{n}(x) \in I_{s} \text{ for some $j$,} \\
    0 & \text{ if } x\in A_{n-1} \setminus \bigcup_{j} U^{1}_{n,j},\\
    t_{n-1}(x)-1 & \text{ if } x\in B_{n-1}\setminus \bigcup_{j} U^{1}_{n,j}.
    \end{cases}
\end{equation*}
Finally let
\[
A_{n}= \{x\in\Delta_{n}: t_{n}(x) = 0\}, \quad B_{n}=
\{x\in\Delta_{n}: t_{n}(x) > 0\}
\]
and
 $$
 A_{n}^\vare= \{x\in\Delta_{n}:
 \dist(f^{n+1}(x),f^{n+1}(A_n))<\vare\}.
 $$
At this point we have completely described the inductive
construction of the sets  $A_n$, $A_n^\vare$,  $B_n$ and
$\{R=n\}$.

We conclude this section with a remark
concerning the role of the sets $B_n$ as a
kind of shield protecting the  sets of the partition  constructed
up to time $n$, and some observations to motivate the last two
sections.

 \subsubsection*{A remark on the construction.}

 Associated to each component $U^0_{n-k}$ of
$\{R=n-k\}$, for some $k>0$, we have a collar $U^1_{n-k}\setminus
U^0_{n-k}$ around it; knowing that the new components of $\{R=n\}$
do not ``intersect too much" $U^1_{n-k}\setminus U^0_{n-k}$ is
important for preventing overlaps on sets of the partition. We
will see that this is indeed the case as long as $\vare>0$ is
taken small enough.

\cle\label{l.claim} If $\vare>0$ is sufficiently small, then
$U_{n}^1\cap\{t_{n-1}>1\}=\emptyset$ for each $U^1_n$. \fle

\dem Take some $k>0$ and let  $U_{n-k}^0$ be a component of
$\{R=n-k\}$ such that its collar $Q_{k}$ (the part of $U_{n-k}^1$
that is mapped by $f^{n-k}$ onto $I_{k}$) intersects $U_{n}^1$.
Recall that $Q_{k}$ is precisely the collar around $U_{n-k}^0$ on
which $t_{n-1}$ takes the value 1. Letting $q_1$ and $q_2$ be any
two points in distinct components of the boundary of $Q_{k}$, we
have
 by Lemma~\ref{l.hyperbolic1} and Remark~\ref{r.c0}
 \begin{equation}\label{e.zq1}
 \dist(f^{n-k}(q_1),f^{n-k}(q_2))\le
 C_0\sigma^{(k-N_0)/2}\dist(f^{n}(q_1),f^{n}(q_2)).
 \end{equation}
We also have
 \begin{eqnarray*}
\dist(f^{n-k}(q_1),f^{n-k}(q_2))&\ge&
\delta_{0}(1+\sigma^{(k-1)/2})-\delta_{0}(1+\sigma^{k/2})\\&=&
\delta_{0} \sigma^{k/2}(\sigma^{-1/2}-1),
 \end{eqnarray*}
which combined with (\ref{e.zq1})  gives
 $$
 \dist(f^{n}(q_1),f^{n}(q_2))\ge C_0^{-1}\sigma^{N_0/2}\delta_{0}(\sigma^{-1/2}-1).
 $$
On the other hand, since $ U^{1}_{n}\subset A_{n-1}^{\varepsilon}$
by construction of $U^{1}_{n}$, taking
 $$
 \vare< C_0^{-1}\sigma^{N_0/2}\delta_{0} (\sigma^{-1/2}-1)
 $$
we  have $U_{n}^1\cap\{t_{n-1}>1\}=\emptyset$.
 \cqd

 \section{The induced map}\label{s.check}

 In this section we briefly discuss the first  and fourth items and prove the second and
 third  items in the statement of the Main Theorem.

 \subsubsection*{The Markov property.}
 The construction detailed in Section \ref{induced map} provides an
 algorithm for the definition of a family of topological balls
 contained in \( \Delta \) and satisfying the Markov property as required.
 In the next two sections we show that this algorithm does indeed
 produce a partition mod 0 of \( \Delta \) and obtain estimates for the
 rate of decay of the tail of the return times.

 \subsubsection*{Uniform expansivity. }

Recall that by construction, the return time \( R(U) \) for \( U
\) an element of the partition \( \mathcal P \) of \( \Delta \),
is formed by a certain number \( n \) of iterations given by the
hyperbolic time of a hyperbolic pre-ball \( V_{n}\supset U \), and
a certain number \( m\leq N_{0} \) of additional iterates which is
the time it takes to go from \( f^{n}(V_{n}) \) which could be
anywhere in \( M \), to \( f^{n+m}(V_{n}) \) which covers \(
\Delta \) completely. By choosing \( R_{0} \) sufficiently large
it then follows from Remark~\ref{r.c0} that there exists a
 constant \( \hat\lambda > 1 \) and a time \( n_{0} \) such that for
 any hyperbolic time \( n\geq n_{0} \) and any point \( x \in V_{n} \)
and $1\le m\le N_0$, we have
\[
\|(Df^{n+m}(x))^{-1}\|^{-1}\geq \hat\lambda > 1
\]
We immediately have the
 uniform expansivity property of the Main Theorem
 \[
 \|(DF_{x})^{-1}\|^{-1}=\|(Df^{R(x)}_{x})^{-1}\|^{-1} \geq
 \hat\lambda > 1.
 \]
In particular, this implies that for any \( x,y \in \Delta \) which
have the same combinatorics, i.e. which remain in the same elements of
the partition \( \mathcal P \) for some number \( s(x,y) \) of
iterates of the induced map \( F \), we have
 \begin{equation}\label{symbolic distance}
     \dist(x,y) \leq \hat\lambda ^{-s(x,y)}.
  \end{equation}

\subsubsection*{Distortion estimates.}
The distortion estimate required for our Main Theorem follows
immediately from \eqref{symbolic distance} above and the following
more classical formulation of the bounded distortion property:

\cle\label{lipdist}
    There exists a constant \( \tilde B > 0 \)
such that for any \( x, y \) belonging to the same element \(
U\in\mathcal P \) with return time \( R \), we have
\[
\log\left|\frac{\det DF(x)}{\det DF(y)}\right|=
\log\left|\frac{\det Df^{R}_{x}}{\det Df^{R}_{y}}\right| \leq \tilde B
\dist(f^{R}(x), f^{R}(y)).
\]
\fle


\begin{proof}
Recall that by
construction, the return time \( R(U) \) for \( U \) an element of the
partition \( \mathcal P \) of \( \Delta \), is formed by a certain
number \( n \) of iterations given by the hyperbolic time of a
hyperbolic pre-ball \( V_{n}\supset U \), and a certain number \(
m\leq N_{0} \) of additional iterates which is the time it takes to go
from \( f^{n}(V_{n}) \) which could be anywhere in \( M \), to \(
f^{n+m}(V_{n}) \) which covers \( \Delta \) completely.
   Some standard formal manipulation
    based on the chain rule gives
    \begin{align*}
\log\left|\frac{\det Df^{R}_{x}}{\det Df^{R}_{y}}\right|
& =
\log\left|\frac{\det Df^{R-n}_{f^n(x)}}{\det Df^{R-n}_{f^n(y)}}\right| +
\log\left|\frac{\det Df^{n}_{x}}{\det Df^{n}_{y}}\right|
    \end{align*}
Since ${f^{i}(x)}$ and ${f^{i}(y)}$ are uniformly bounded away
from $\cs$ for $n\le i\le R$ (recall Remark~\ref{r.c0}), we may
write
$$
\log\left|\frac{\det Df^{R-n}_{f^n(x)}}{\det Df^{R-n}_{f^n(y)}}\right|
 \le
B_1 \dist(f^{R}(x),f^{R}(y))
$$
where $B_1$ is some constant not depending on $x,y$ or $R$.
On the other hand, by construction of $V_n$ (see the proof of Lemma 5.2
in \cite{ABV}), there must be some $z\in V_n$ for which $n$ is a
hyperbolic time and such that, for $0\le j <n$, the distance from
$f^j(z)$ to either $f^j(x)$ or $f^j(y)$ is smaller than
 $\dist(f^{n}(x),f^{n}(y)) \sigma^{(n-j)/2}$, which is much smaller than
$\sigma^{b(n-j)}\le \dist(f^j(z),\cs)$. Thus, by (S3) we have
 $$
\log\left|\frac{\det Df^{n}_{x}}{\det Df^{n}_{y}}\right|
\le
\sum_{j=0}^{n-1}\log\left|\frac{\det Df_{f^j(x)}}{\det Df_{f^j(y)}}\right|
\le
\dist(f^{n}(x),f^{n}(y)) \sum_{j=0}^{n-1}2B
\frac{\sigma^{(n-j)/2}}{\sigma^{b\beta(n-j)}}.
 $$
Since $b\beta<1/2$, there must be some $B_2>0$ such that
 $$
\log\left|\frac{\det Df^{n}_{x}}{\det Df^{n}_{y}}\right|
\le
 B_2\dist(f^{n}(x),f^{n}(y))
 .$$
Using again that ${f^{i}(y)}$ and ${f^{i}(y)}$ are uniformly
bounded away from $\cs$ (for $n\le i\le R$ (cf. Remark~\ref{r.c0})
it follows that
 $$
 \dist(f^{n}(x),f^{n}(y))
\le B_2 \dist(f^{R}(x),f^{R}(y)),
 $$
 where $B_2$ is some constant not depending on $x,y$ or $R$.
This completes the proof of the lemma.
\end{proof}

\subsubsection*{Looking ahead: probabilistic estimates.}

 For proving the Main Theorem we only need to study the
decay of $\leb(\Delta_n)$ in terms of $\leb(\Gamma_n)$. That is
our purpose in the next two sections. We will show in
Proposition~\ref{p.ab} that there is a constant $a_0>0$ such that
 \begin{equation}\label{e.ab}
  \leb(B_{n})\leq a_0 \leb(A_{n}).
 \end{equation}
We will also show in Proposition~\ref{p.construction} that there
are $N=N(\vare)\ge 1$ and a  constant $c_0>0$ for which
  \begin{equation}\label{e.construction}
\leb\left(\bigcup_{i=0}^N\big\{R=n+i\big\}\right)\ge c_0
\leb(A_{n-1}\cap H_{n}).
 \end{equation}
Taking into account that $\Delta_n=A_n\cup B_n$, it easily follows
from (\ref{e.ab}) and (\ref{e.construction}) that there is a
constant $b_0>0$ such that
 $$
 \leb\left(\bigcup_{i=0}^N\big\{R=n+i\big\}\right)\ge b_0
 \frac{\leb(A_{n-1}\cap H_{n})}{\leb(A_{n-1})}\leb(\Delta_{n-1}).
 $$
This immediately implies that
 \begin{equation}\label{eq.recorre}
 \leb\left(\Delta_{n+N}\right)\le \left(1-b_0
 \frac{\leb(A_{n-1}\cap H_{n})}{\leb(A_{n-1})}\right)\leb(\Delta_{n-1}).
 \end{equation}
It is no restriction to assume that $R_0>2(N+1)$ and we do it.
Take any large $n$ and let $k_0\ge1$ be the smallest integer for
which $n-1-k_0(N+1)\le R_0$. The above assumption on $R_0$ and $N$
implies that $ n-(k_0+1)(N+1)\ge 1$. Now we consider the partition
of $\{n-(k_0+1)(N+1),\dots,n-1\}$ into the sets
  \begin{eqnarray*}
    J_N &= &\{n-1,n-1-(N+1),\dots,n-1-k_0(N+1)\},\\
     &\vdots& \\
  J_1 &= &\{n-N,n-N-(N+1)\dots,n-N-k_0(N+1)\},\\
 J_{0} &= &\{ n-(N+1),n-2(N+1)\dots,n-(k_0+1)(N+1)\}.
 \end{eqnarray*}
 Applying (\ref{eq.recorre})  repeatedly we arrive at the following set
 of $N+1$ inequations:
  \begin{eqnarray*}
  \leb\left(\Delta_{n+N}\right)
   &\le &
  \prod_{j\in J_N}\left(1-b_0 \frac{\leb(A_{j}\cap
  H_{j+1})}{\leb(A_{j})}\right)\leb(\Delta_{0}),\\
   \vdots\quad\quad && \\
   \leb\left(\Delta_{n}\right)
   &\le &
  \prod_{j\in J_0}\left(1-b_0 \frac{\leb(A_{j}\cap
  H_{j+1})}{\leb(A_{j})}\right)\leb(\Delta_{0}).
  \end{eqnarray*}
Multiplying the terms in the inequations above and ignoring
factors from $n-(k_0+1)(N+1)$ to $R_0-1$ on the right hand side
(observe that those factors are smaller than 1), we obtain
 \begin{equation}\nonumber
 \prod_{j=0}^{N}\leb\left(\Delta_{n+j}\right)\le \prod_{j=R_0}^{n-1} \left(1-b_0
 \frac{\leb(A_{j-1}\cap
 H_{j})}{\leb(A_{j-1})}\right)\leb(\Delta_{0})^{N+1}.
 \end{equation}
Taking into account that $(\Delta_n)_n$ forms a decreasing
sequence of sets we finally have
 \begin{equation}\label{e.media}
 \leb\left(\Delta_{n+N}\right)\le \exp\left(-\frac{b_0}{N+1}\sum_{j=R_0}^{n}
 \frac{\leb(A_{j-1}\cap
 H_{j})}{\leb(A_{j-1})}\right)\leb(\Delta_{0}).
 \end{equation}
 In Section~\ref{s.decay} we
will prove the Main Theorem  by considering several
different cases, according to the behavior of the proportions
${\leb(A_{n-1}\cap H_{n})}/{\leb(A_{n-1})}$.
 It is not hard to check that if the average
 $$\frac1n\sum_{j=1}^n\frac{\leb(A_{j-1}\cap H_{j})}{\leb(A_{j-1})}$$
 is bounded away from 0 for
large $n$, then  $\leb(\Delta_n)$ decays exponentially fast to 0.
This happens, for instance,  when $f$ is uniformly expanding.

\section{Transitional metric estimates}
\label{transitional estimates}

The goal of this section is to prove several estimates
relating the Lebesgue measure of the sets $A_n$, $A_n^\vare$, $B_n$
and $\{R=n\}$.  The first result shows that a fixed proportion of
$A_{n-1}$ having $n$ as a hyperbolic time gives rise to new elements
of the partition with return time not exceeding $n$ too much.  We
discuss the relative proportion of the sets \( A_{n} \) and \( B_{n}
\) in \( \Delta_{n} \).


    \cpr\label{p.construction}
    There exist $c_0>0$ and $N=N(\vare)$ such that for every \( n\ge1\)
    \[
     \leb\left(\bigcup_{i=0}^N\big\{R=n+i\big\}\right)\ge c_0 \leb(A_{n-1}\cap H_{n}).
    \]
    \fpr
\begin{proof} Take $r=5\delta_0C_0^{N_0}$, where $N_0$ and
$C_0$ are given by Lemma~\ref{le:grow2} and Remark~\ref{r.c0},
respectively. Let $\{z_j\}$ be a maximal set in $f^n(A_{n-1}\cap
H_n)$ with the property that $B(z_j,r)$ are pairwise disjoint. By
maximality we have
 $$
 \bigcup_j B(z_j,2r)\supset f^n(A_{n-1} \cap H_n).
 $$
Let $x_j$ be a point in $H_n$ such that $f^n(x_j)=z_j$ and
consider the hyperbolic pre-ball $V_n(x_j)$ associated to $x_j$.
Observe that $f^n$ sends $ V_n(x_j)$ diffeomorphically onto a ball
of radius $\delta_1$ around $z_j$ as in Lemma~\ref{l.hyperbolic1}.
In what follows, given $B\subset B(z_j,\delta_1)$, we will simply
denote
 $(f^{n}\vert V_n(x_j))^{-1}(B)$ by $f^{-n}(B)$.

Our aim now is to prove that $f^{-n}(B(z_j, r))$ contains some
component of $\{R=n+k_j\}$  with  $0\le k_j\le N_\vare+N_0$. We
start by showing that
 \begin{equation}\label{eq.tnk}
  t_{n+k_j}\vert f^{-n}(B(z_j,\vare))>0\quad\text{for some $0\le k_j\le
  N_\vare+N_0$}.
  \end{equation}
Assume by contradiction that
  $
  t_{n+k_j}\vert f^{-n}(B(z_j,\vare))=0$ for all $0\le k_j\le
  N_\vare+N_0
  $.
 This implies that $f^{-n}(B(z_j,\vare))\subset A_{n+k_j}^\vare$
 for all $0\le k_j\le N_\vare+N_0$.
 Using Lemma~\ref{le:grow1} we may find  a hyperbolic pre-ball
 \( V_{m}\subset B(z_j,\vare) \) with \( m\leq  N_{\varepsilon} \).
 Now, since $f^m(V_{m})$ is a ball $B$ of radius $\delta_1$ it follows
from Lemma~\ref{le:grow2} that there is some $V\subset B$
and $m'\le N_0$ with $f^{m'}(V)=\Delta_0$. Thus, taking $k_j=m+m'$
we have that $0\le k_j\le N_\vare+N_0$ and $f^{-n}(V_m)$ is an
element of $\{R=n+k_j\}$ inside $f^{-n}(B(z_j,\vare))$. This
contradicts the fact that $t_{n+k_j}\vert f^{-n}(B(z_j,\vare))=0$
 for all $0\le k_j\le N_\vare+N_0$, and so (\ref{eq.tnk}) holds.

Let $k_j$ be the
 smallest integer $0\le k_j\le N_\vare+N_0$ for which
$t_{n+k_j}\vert f^{-n}(B(z_j,\vare))>0$.
Since
  $$
  f^{-n}(B(z_j,\vare))\subset A_{n-1}^\vare\subset \{t_{n-1}\le1 \},
  $$
 there must be some element
$U^{0}_{n+k_j}(j)$ of $\{R=n+k_j\}$ for which
 $$
 f^{-n}(B(z_j,\vare))\cap U_{n+k_j}^1(j)\neq\emptyset.
 $$
Recall that by definition $f^{n+k_j}$ sends $U_{n+k_j}^1(j)$
diffeomorphically onto $\Delta_0^1$, the ball  of radius
$(1+s)\delta_0$ around $p$. From time $n$ to $n+k_j$ we may have
some final ``bad" period of length at most $N_0$ where the
derivative of $f$ may contract, however being bounded from below
by $1/C_0$ in each step. Thus, the diameter of
$f^n(U_{n+k_j}^1(j))$ is at most $4\delta_0C_0^{N_0}$. Since
$B(z_j,\vare)$ intersects $f^n(U_{n+k_j}^1(j))$ and
$\vare<\delta_0<\delta_0C_0^{N_0}$, we have by definition of $r$
 $$
 f^{-n}(B(z_j,r))\supset U_{n+k_j}^0(j).
 $$
Thus we have shown  that $f^{-n}(B(z_j, r))$ contains some
component of $\{R=n+k_j\}$ with  $0\le k_j\le N_\vare+N_0$.
Moreover, since $n$ is a hyperbolic time for $x_j$, we have by the
distortion control given by Lemma~\ref{l.hyperbolic1}
 \begin{equation}\label{eq.quo1}
 \frac{ \leb(f^{-n}(B(z_j,2r)))}{\leb(f^{-n}(B(z_j,r)))}
 \le
 {D_1}\frac{\leb(B(z_j,2r))}{\leb(B(z_j,r))}
 \end{equation}
and
 \begin{equation}
\label{eq.quo2} \frac{\leb(f^{-n}(B(z_j,r)))}{\leb(
U_{n+k_j}^0(j))}
 \le
 {D_1}\frac{\leb(B(z_j,r))}{\leb(f^n(U_{n+k_j}^0(j)))}.
 \end{equation}
Here we are implicitly assuming that $2r<\delta_1$. This can be
done just by taking $\delta_0$ small enough. Note that the
estimates on $N_0$ and $C_0$ improve when we diminish $\delta_0$.

From time $n$ to time $n+k_j$ we have at most $k_j=m_1+m_2$
iterates with $m_1\le N_\vare$, $m_2\le N_0$ and
$f^n(U_{n+k_j}^0(j)))$ containing some point $w_j\in H_{m_1}$. By
the definition of $(\sigma,\delta)$-hyperbolic time we have that
$\dist_\delta (f^i(x),\cs)\ge \sigma^{bN_\vare}$ for every $0\le
i\le m_1$, which by the uniform distortion control implies that
there is some constant $D=D(\vare)>0$ such that $|\det
(Df^i(x))|\le D$ for  $0\le i\le m_1$ and $x\in
f^n(U_{n+k_j}^0(j))$. On the other hand, since the first $N_0$
preimages of $\Delta_0$ are uniformly bounded away from $\cs$ we
also have some $D'>0$ such that $|\det (Df^i(x))|\le D'$ for every
$0\le i\le m_2$ and $x$ belonging to an $i$ preimage of
$\Delta_0$. Hence,
$$\leb(f^n(U_{n+k_j}^0(j)))\ge \frac{1}{DD'}\leb(\Delta_0),$$
which combined with (\ref{eq.quo2}) gives
 $$
 \leb(f^{-n}(B(z_j,r)))\le C\leb(
U_{n+k_j}^0(j)),
 $$
with $C$ only depending on $D_1$, $D$, $D'$, $\delta_0$ and the
dimension of $M$. We also deduce from (\ref{eq.quo1}) that
 $$
 \leb(f^{-n}(B(z_j,2r)))\le C'\leb(
f^{-n}(B(z_j,r)))
 $$
with $C'$ only depending on $D_1$ and the dimension of $M$.
 Finally let us compare the Lebesgue measure of
the sets $\bigcup_{i=0}^N\big\{R=n+i\big\}$ and $A_{n-1}\cap
H_{n}$. We have
 $$
 \leb\big(A_{n-1}\cap
 H_n \big)\le \sum_j \leb(f^{-n}(B(z_j,2r)))
 \le C'
 \sum_j \leb(f^{-n}(B(z_j,r))).
  $$
On the other hand, by the disjointness of the balls $B(z_j,r)$ we
have
$$
 \sum_j \leb(f^{-n}(B(z_j,r)))\le C
 \sum_j \leb( U_{n+k_j}^0(j)) \le C
 \leb\left(\bigcup_{i=0}^N\big\{R=n+i\big\}\right).
  $$
We just have to take $c_0^{-1}=CC'$.
\end{proof}

Let us prove now a couple of useful lemmas. The first one gives a
lower bound for the flow of mass from $B_{n-1}$ to $A_n$, and
second one gives a lower bound for the flow of mass from $A_{n-1}$
to $B_n$ and $\{R=n\}$.

\cle\label{l.flowb} There exists $a_1>0$  such that for every
$n\ge1$
$$\leb(B_{n-1}\cap A_n)\ge a_1\leb(B_{n-1}).$$
Moreover, $a_1$ is bounded away from 0 independently from
$\delta_0$.
 \fle
  \dem
  It is enough to see this for each component of $B_{n-1}$
 at a time. Let  $C$ be a component of $B_{n-1}$ and $Q$
be  its outer ring corresponding to $t_{n-1}=1$. Observe that by
Lemma~\ref{l.claim} we have $Q=C\cap A_n$. Moreover, there must be
some $k<n$ and a component $U^0_k$ of $\{R=k\}$ such that $f^k$
maps $C$ diffeomorphically onto $\bigcup_{i=k}^\infty I_i$ and $Q$
onto $I_k$, both with uniform bounded distortion (not depending on
$\delta_0$ or $n$).
 Thus, it is sufficient to compare the Lebesgue measures of $\bigcup_{i=k}^\infty
I_i$ and $ I_k$. We have
 $$
 \frac{\leb( I_k)}{\leb(\bigcup_{i=k}^\infty
I_i)}\thickapprox\frac{
[\delta_0(1+\sigma^{(k-1)/2})]^d-[\delta_0(1+\sigma^{k/2})]^d}
{[\delta_0(1+\sigma^{(k-1)/2})]^d-\delta_0^d}\thickapprox
1-\sigma^{1/2}.
 $$
 Clearly this proportion does not depend on  $\delta_0$.
  \cqd

The second item of the lemma below is apparently counterintuitive,
since our main goal is to make the points in $\Delta_0$ have small
return times. However, this is needed for keeping $\leb(A_n)$
uniformly much bigger than $\leb(B_n)$. This will help us in the
statistical estimates of the last section.

\cle\label{l.flowa} There exist $b_1=b_1(\delta_0)>0$ and
$c_1=c_1(\delta_0)>0$ with $b_1+c_1<1$ such that for every $n\ge1$
\begin{enumerate}
\item $\leb(A_{n-1}\cap B_n)\le b_1\leb(A_{n-1})$;
\item $\leb(A_{n-1}\cap \{R=n\})\le c_1\leb(A_{n-1})$.
\end{enumerate}
    Moreover, $b_1\to 0$ and $c_1\to 0$ as $\delta_0\to 0$.
\fle \dem It is enough to prove this for each  neighborhood of a
component $U^0_n$ of $\{R=n\}$. Observe that by construction we
have $U^3_n\subset A_{n-1}^\vare$, which means that $U^2_n \subset
A_{n-1}$, because $\vare<\delta_0<\sqrt\delta_0$. Using the
uniform bounded distortion of $f^n$ on $U^3_n$ given by
Lemma~\ref{le:grow1} and Lemma~\ref{le:grow2} (cf.
Remark~\ref{r.c0}) we obtain
 $$
 \frac{\leb(U^1_n\setminus U^0_n)}{\leb(U^2_n\setminus U^1_n)}
 \thickapprox
 \frac{\leb(\Delta^1_0\setminus \Delta^0_0)}{\leb(\Delta^2_0\setminus \Delta^1_0)}
 \thickapprox
 \frac{\delta_0^d}{\delta_0^{d/2}}\ll 1,
 $$
which gives the first estimate. Moreover,
$$
 \frac{\leb( U^0_n)}{\leb(U^2_n\setminus U^1_n)}
 \thickapprox
 \frac{\leb( \Delta^0_0)}{\leb(\Delta^2_0\setminus \Delta^1_0)}
 \thickapprox
 \frac{\delta_0^d}{\delta_0^{d/2}}\ll 1,
 $$
and this gives the second one.
 \cqd

The next result is a consequence of the estimates we obtained in
the last two lemmas. The proof is essentially the same of the
uniformly hyperbolic case; see \cite{Y1}. Here we need to be more
careful on the estimates.

\cpr\label{p.ab}
    There exists  $a_0=a_0(\delta_0)>0$ such that for every $n\ge1$
    \[
    \leb(B_{n}) \leq a_0\leb(A_{n}) .
    \]
    Moreover, $a_0\to 0$ as $\delta_0\to 0$.
\fpr
\begin{proof}
We have  by Lemma~\ref{l.flowa}
 \begin{equation}\label{eq.eta}
 \leb(A_n)\ge
 (1-b_1-c_1)\leb(A_{n-1}).
 \end{equation}
Letting $\eta=1-b_1-c_1$ we define
 $$\widehat a=\frac {b_1+c_1}{a_1}\qand
 a_0=\frac{(1+a_1)b_1+c_1}{a_1\eta}.
 $$
 The fact that $a_0\to 0$ when $\delta_0\to 0$ is a consequence of
 $b_1\to 0$ and $c_1\to 0$ when $\delta_0\to 0$ and $a_1$ being bounded
 away from 0.
 Observe that $0<\eta<1$ and $\widehat a< a_0$.
Now  the proof of  the proposition follows by induction. The
result obviously holds for $n$ up to $R_0$. Assuming that it holds
for $n-1\ge R_0$ we will show that it also holds for $n$, by
considering separately the cases $\leb(B_{n-1})>\widehat a
\leb(A_{n-1})$ and $\leb(B_{n-1})\le\widehat a \leb(A_{n-1}).$

Assume first that  $\leb(B_{n-1})>\widehat a \leb(A_{n-1})$.
 We may write
  $$\leb(B_{n-1})=\leb(B_{n-1}\cap A_n)+\leb(B_{n-1}\cap B_n),
  $$
which by Lemma~\ref{l.flowb} gives
\begin{equation}\label{e.bb}
  \leb(B_{n-1}\cap B_n)\le (1-a_1)\leb(B_{n-1}).
\end{equation}
Since we also have
 $$\leb(B_{n})=\leb(B_{n}\cap B_{n-1})+\leb(B_{n}\cap A_{n-1}),
 $$
 it follows from (\ref{e.bb}) and Lemma~\ref{l.flowa} that
  $$
  \leb(B_n)\le (1-a_1)\leb(B_{n-1})+b_1\leb(A_{n-1}) ,
  $$
  which according to the case we are considering leads to
  \begin{eqnarray*}
  \leb(B_n)\le (1-a_1)\leb(B_{n-1})+\frac{b_1a_1}{b_1+c_1}\leb(B_{n-1})
  < \leb(B_{n-1}).
  \end{eqnarray*}
On the other hand, we have by (\ref{eq.eta}) that
$\leb(A_n)<\leb(A_{n-1})$, which together with the last inequality
and the inductive hypothesis yields
 $$\frac{\leb(B_n)}{\leb(A_n)}<\frac{\leb(B_{n-1})}{\leb(A_{n-1})}\le a_0,$$
 which gives the result in the first case.

Assume now that $\leb(B_{n-1})\le\widehat a \leb(A_{n-1}).$
 Since we have
 $$\leb(B_{n})=\leb(B_{n}\cap B_{n-1})+\leb(B_{n}\cap A_{n-1}),
 $$
 it follows from Lemma~\ref{l.flowa} that
  $$\leb(B_{n})\le\leb( B_{n-1})+b_1\leb( A_{n-1}).
 $$
Hence
$$\frac{\leb(B_n)}{\leb(A_n)}<
\frac{\leb( B_{n-1})+b_1\leb( A_{n-1})}{\eta\leb(A_{n-1})} \le
\frac{\widehat a+b_1}{\eta}=a_0,
$$
which gives the result also in this case.
\end{proof}

It will be useful to establish the following  consequence of the
last two results.

\cco\label{c.omega}
    There exists  $c_2>0$ such that for every $n\ge1$
    \[
    \leb(\Delta_{n}) \le c_2\leb(\Delta_{n+1}) .
    \]
\fco
 \dem
By Lemma~\ref{l.flowa} we have
 $$
 \leb(\Delta_{n+1})\ge \leb(A_{n+1}) \ge (1-b_1-c_1)\leb(A_n).
 $$
On the other hand, by Proposition~\ref{p.ab},
 $$\leb(\Delta_n)=\leb(A_n)+\leb(B_n)\le
 (1+a_0^{-1})\leb(A_n).
 $$
It is enough to take $c_2=(1+a_0^{-1})/(1-b_1-c_1)$.
  \cqd

\section{Asymptotic metric estimates}
\label{s.decay}

We start this section by recalling that  $\theta>0$ was obtained
in Lemma~\ref{l:hyperbolic2} and gives a lower bound for the
frequency of hyperbolic times; it  only depends on the non-uniform
expansion coefficient $\lambda$ and the map $f$.

Before we go into the main proposition of this section which will
enable us to conclude the proof of the Main Theorem, let us impose
one more requirement on the choice of $\delta_0$: let $\gamma>0$
be some positive number (to be determined later) and take $
0<\alpha<(\theta/{12})^{\gamma+1}$. Then we choose $\delta_0>0$
small so that $a_0=a(\delta_0)<2\alpha$.

We define for each $n\ge 1$
 $$
 E_n=\left\{ j\le n\colon \frac{\leb(A_{j-1}\cap
 H_j)}{\leb(A_{j-1})}<\alpha\right\},
 $$
and
 $$F=\left\{n\in\NN\colon \frac{\#
 E_n}{n}>1-\frac{\theta}{12}\right\}.$$

\cpr\label{p.final} Take any  $n\in F$ with $n\ge R_0> 12/\theta$.
If $\leb(A_n)\ge 2 \leb(\Gamma_n)$, then there is some
$0<k=k(n)<n$ for which
 $$\frac{\leb(A_n)}{\leb(A_k)}<\left(\frac kn\right)^\gamma.
 $$
 \fpr

 \dem
We have for $j\le n$
 $$
 \frac{\leb(A_{n}\cap H_j)}{\leb(A_{n})}\ge
 \frac{\leb(A_{n}\setminus \Gamma_n)}{\leb(A_n)}\cdot
 \frac{\leb((A_{n}\setminus \Gamma_n)\cap H_j)}{\leb(A_{n}\setminus \Gamma_n)}
 \ge \frac12\cdot\frac{\leb((A_{n}\setminus \Gamma_n)\cap
 H_j)}{\leb(A_{n}\setminus \Gamma_n)},
 $$
which together with the conclusion of
Corollary~\ref{c:hyperbolic3} for the set $A_{n}\setminus
\Gamma_n$ gives
 \begin{equation}\label{eq.media}
 \frac1n\sum_{j=1}^n\frac{\leb(A_{n}\cap
 H_j)}{\leb(A_{n})} \ge \frac\theta2.
 \end{equation}
Let
 $$G_n=\left\{j\in E_n\colon
 \frac{\leb(A_{j-1})}{\leb(A_n)}>\frac{\theta}{12\alpha}\right\}.
 $$
Since $n\in F$, we have
 \begin{eqnarray*}
 \frac1n\sum_{j=1}^n\frac{\leb(A_{n}\cap
 H_j)}{\leb(A_{n})}
 &\le&
 \frac\theta{12}+\frac1n\sum_{j\in E_n}\frac{\leb(A_{n}\cap
 H_j)}{\leb(A_{n})}\\
 &\le&
 \frac\theta{12}+\frac1n\sum_{j\in E_n\setminus G_n}\frac{\leb(A_{n}\cap
 H_j)}{\leb(A_{n})}+\frac{\#G_n}n. \\
 \end{eqnarray*}
Now, for $j\in E_n\setminus G_n$,
 \begin{eqnarray*}
 \frac{\leb(A_{n}\cap H_j)}{\leb(A_{n})}
 &=&
 \frac{\leb(A_{n}\cap H_j)}{\leb(A_{j-1})}\cdot\frac{\leb(A_{j-1})}{\leb(A_{n})}\\
 &\le &
 \left(\frac{\leb(A_{n}\cap A_{j-1}\cap H_j)}{\leb(A_{j-1})}+
 \frac{\leb((A_{n}\setminus A_{j-1})\cap H_j)}{\leb(A_{j-1})}\right)
 \frac{\leb(A_{j-1})}{\leb(A_{n})}\\
 &\le &
 \left(\frac{\leb( A_{j-1}\cap H_j)}{\leb(A_{j-1})}+
 a_0\right)
 \frac{\theta}{12\alpha}.
 \end{eqnarray*}
For this last inequality we used the fact that $(A_{n}\setminus
A_{j-1})\subset B_{j-1}$ and  $j\notin G_n$. Hence
 \begin{eqnarray*}
 \frac1n\sum_{j=1}^n\frac{\leb(A_{n}\cap H_j)}{\leb(A_{n})}
 &\le &
 \frac\theta{12}+ \frac1n\sum_{j\in E_n\setminus G_n}\frac{\leb( A_{j-1}\cap
 H_j)}{\leb(A_{j-1})}\frac{\theta}{12\alpha}+ a_0
 \frac{\theta}{12\alpha}+\frac{\#G_n}{n}\\
 &<&
 \frac\theta{12}+ \alpha\frac{\theta}{12\alpha}+ a_0
 \frac{\theta}{12\alpha}+\frac{\#G_n}n.
 \end{eqnarray*}
 By the choice of $a_0$  we  have that the third term
 in the last sum above is smaller than
 $\theta/6$. So, using (\ref{eq.media}) we obtain
 \begin{equation}\label{fr.v}
 \frac{\#G_n}n> \frac\theta6.
 \end{equation}
Now, defining
 $$
 k=\max(G_n)-1,
 $$
 we have
\begin{eqnarray*}
\leb(A_n)<\frac{12\alpha}\theta\leb(A_k).
 \end{eqnarray*}
 It follows from (\ref{fr.v}) that
 $ k+1> \theta n/6 $,
 and so
 $ k/n>\theta/{12}$,
 because $n\ge R_0> 12/\theta$.
 Since we have chosen
 $\alpha<\left(\theta/12\right)^{\gamma+1}, $
it follows  that
 $$
 \left(\frac kn\right)^{\gamma}>\frac{12}{\theta}\left(\frac\theta{12}\right)^{\gamma+1}
 >\frac{12\alpha}\theta.
 $$
This completes the proof of the result. \cqd

Now we are ready to conclude the proof of our Main Theorem, namely
the decay estimate on the tail of return times. Observe that by
Proposition~\ref{p.ab} we have $\leb(\Delta_n)\lesssim \leb(A_n)$,
and so it is enough to derive the tail estimate for $\leb(A_n)$ in
the place of $\leb(\{R>n\})=\leb(\Delta_n)$.

Given any large integer $n$, we consider  the following
situations:
\begin{enumerate}
\item If $n\in\NN\setminus F$, then by (\ref{e.media}) and
Corollary~\ref{c.omega} we have
 \begin{equation}\label{e.f}\nonumber
 \leb(\Delta_n)\le c_2^N\exp\left({-\frac{b_0\theta\alpha}{12(N+1)}
 (n-R_0)}\right){\leb(\Delta_0)}.
 \end{equation}
\item If $n\in F$, then we distinguish the following two
cases:
\begin{enumerate}
\item If $\leb(A_n)< 2 \leb(\Gamma_n)$, then nothing has to be
done.
\item If $\leb(A_n)\ge 2 \leb(\Gamma_n)$, then we apply
Proposition~\ref{p.final} and get some $k_1<n$ for which
$${\leb(A_n)}<\left(\frac {k_1}n\right)^\gamma {\leb(A_{k_1})}.
 $$
\end{enumerate}
\end{enumerate}
The only situation we are left to consider is 2(b). In such case,
either $k_1$ is in situation~1 or 2(a), or by
Proposition~\ref{p.final} we can find $k_2<k_1$ for which
 $$
 {\leb(A_{k_1})}<\left(\frac {k_2}{k_1}\right)^\gamma {\leb(A_{k_2})}.
 $$
Arguing inductively we show that there is a sequence of integers
$0<k_s<\cdots<k_1<n$ for which   one of the following cases
eventually  holds.
 \begin{enumerate}
\item[(I)]\quad $\displaystyle{\leb(A_{n})<\left(\frac {k_s}{n}\right)^\gamma
c_2^N\exp\left({-\frac{b_0\theta\alpha}{12(N+1)}
 (k_s-R_0)}\right){\leb(\Delta_0)}}$.
\item[(II)]\quad $\displaystyle{\leb(A_{n})<\left(\frac {k_s}{n}\right)^\gamma
\leb(\Gamma_{k_s})}$.
\item[(III)]\quad $\displaystyle{\leb(A_{n})<\left(\frac {R_0}{n}\right)^\gamma
\leb(\Delta_0)}$.
 \end{enumerate}
Case (III) corresponds to falling into situation 2(b) repeatedly
until $k_s\le R_0$. Observe that until now $\gamma>0$ is
arbitrary. So, the only case we are left to consider is (II).

Assume  that $\leb(\Gamma_n)\le\mathcal O(n^{-\gamma})$ for some
$\gamma>0$. In this case there must be some $C>0$ such that
 $
 {k^\gamma}\leb(\Gamma_k)\le {C}$ for all $k\in \NN$,
 which applied to $k_s$ in case (II) leads to
  $\leb(A_n)\le\mathcal O(n^{-\gamma})$.


\end{document}